\documentclass[10pt]{article}
\usepackage{cite}
\usepackage{enumitem}
 \usepackage{amsfonts}
\usepackage{amsmath}
\usepackage{amsthm}
\usepackage{amssymb}
\usepackage{amscd}
\newtheorem{definition}{Definition}[section]
\newtheorem{theorem}{Theorem}[section]

\newtheorem{lemma}{Lemma}[section]

\newtheorem{remark}{Remark}[section]
\newtheorem{example}{Example}[section]

\begin{document}
\title{{Lower semicontinuity of solution mappings for parametric fixed point problems with  applications}\thanks{This work was supported by  the National Natural Science Foundation of China (11471230, 11671282).}}
\author{{Yu Han$^{a}$ and  Nan-jing Huang$^b$ \footnote{Corresponding author,  E-mail: nanjinghuang@hotmail.com; njhuang@scu.edu.cn} }\\
{\small\it a. Department of Mathematics, Nanchang University, Nanchang,
Jiangxi 330031, P.R. China}\\
{\small\it b. Department of Mathematics, Sichuan University, Chengdu,
Sichuan 610064, P.R. China}}
\date{}
\maketitle
\vspace*{-9mm}
\begin{center}
\begin{minipage}{5.8in}
{\bf Abstract.} In this paper, we establish the lower semicontinuity of the solution mapping and of the approximate solution mapping for parametric
 fixed point problems under some suitable conditions. As applications, the lower semicontinuity result applies to the parametric vector quasi-equilibrium problem, and allows to prove the existence of solutions for generalized Stackelberg games.

\end{minipage}
\end{center}
\section{Introduction}
The  semicontinuity of solution mappings of vector equilibrium problems has been investigated by several authors, see \cite{Anh1, Anh2, Anh3, Anh4, ChenCR, GongU, HG, HH, HH2, Huang, Li1} and the references therein.  Recently, in order to show the semicontinuity of the solution mappings for the parametric (vector) quasi-equilibrium problems, all the solution mappings of the parametric fixed point problems are assumed to be lower semicontinuous in the literature \cite{Anh1, Anh2, Anh3}.  We note that in the literature mentioned above, the authors have not given any conditions to guarantee the lower semicontinuity of the solution mappings of the parametric fixed point problems. On the other hand,  it is difficult to obtain the explicit solutions for some real problems when the data concerned with the problems are perturbed by noise and so the mathematical models are usually solved by numerical methods for approximating the exact solutions. Therefore, one nature question is:  can we provide conditions ensuring the lower semicontinuity of the (approximate) solution mappings?

The main purpose of this paper is to make a new attempt to establish the lower semicontinuity of the solution mapping and of the approximate solution mapping for parametric fixed point problems under suitable conditions.  The rest of the paper is organized as follows. Section 2 presents some necessary notations and lemmas. In Section 3, we establish the lower semicontinuity of the solution mapping and of the approximate solution mapping for parametric fixed point problems.  In Section 4, the lower semicontinuity result applies to the parametric vector quasi-equilibrium problem, and allows to prove the existence of solutions for generalized Stackelberg games.

\section{Preliminaries}
 Throughout this paper, unless otherwise specified,   let $\Lambda$ and  $X$   be two normed vector spaces,  $ {\mathbb{R}_ + } = \left\{ {x \in \mathbb{R}:x \ge 0} \right\}$, $ {\mathbb{R}_ + ^0 } = \left\{ {x \in \mathbb{R}:x > 0} \right\}$ and  $\mathbb{N}=\{1,2,\cdots\}$.
Let $A$ be a nonempty subset of $X$ and $T:A \times \Lambda  \to {2^A}$ be a set-valued mapping. For $\lambda  \in \Lambda $, we consider the following  parametric fixed point problem  consisting of finding ${x_0} \in A$ such that
\begin{center}(PFPP) \quad  ${x_0} \in T\left( {{x_0},\lambda } \right)$. \end{center}

For $\lambda  \in \Lambda $, let ${S}\left( { \lambda } \right)$ denote the set of all  solutions of (PFPP), i.e.
$$S\left( \lambda  \right) = \left\{ {x \in A:x \in T\left( {x,\lambda } \right)} \right\}.$$

For $\left( {\lambda ,\varepsilon } \right) \in \Lambda  \times {\mathbb{R}_ + }$, let $E\left( {\lambda ,\varepsilon } \right)$ denote the set of all  $\varepsilon$-approximate   solutions of (PFPP), i.e.
$$E\left( {\lambda ,\varepsilon } \right) = \left\{ {x \in A:d\left( {x,T\left( {x,\lambda } \right)} \right) \le \varepsilon } \right\},$$
where   $d\left( {x,T\left( {x,\lambda } \right)} \right) = \mathop {\inf }\limits_{y \in T\left( {x,\lambda } \right)} d\left( {x,y} \right)$ and $d\left( {x,y} \right) = \left\| {x - y} \right\|$.

Denote  the boundary of $D$ by ${\partial D}$,  the complement of $D$ by ${D^c}$, the closure of $D$ by ${\rm{cl}}D$ and the interior of $D$ by ${\mathop{\rm int}} D$.

\begin{definition}  (\cite{Holmes})  A nonempty  convex  subset $D$ of $X$ is said to be rotund if the boundary of $D$ does not contain line segments, i.e., for any ${x_1},{x_2} \in D$ with ${x_1} \ne {x_2}$, $\left( {{x_1},{x_2}} \right) \cap {\left( {\partial D} \right)^c} \ne \emptyset $, where $\left( {{x_1},{x_2}} \right) = \left\{ {\lambda {x_1} + \left( {1 - \lambda } \right){x_2}:\lambda  \in \left( {0,1} \right)} \right\}.$
\end{definition}
\begin{remark} Let $D$ be  a nonempty  convex  subset  of $X$. Then it is easy to see  that $D$  is  rotund if and only if for any ${x_1},{x_2} \in D$ with ${x_1} \ne {x_2}$, there exists ${\lambda _0} \in \left( {0,1} \right)$ such that ${\lambda _0}{x_1} + \left( {1 - {\lambda _0}} \right){x_2} \in {\mathop{\rm int}} D$.
 Let $D = \left\{ {\left( {x,y} \right) \in {{\mathbb{R}}^2}:{x^2} + {y^2} \le 1} \right\}$. Then it is clear that $D$  is  rotund.
\end{remark}

\begin{definition} Let $\Delta$ and $\Delta_1$ be two topological vector spaces. A set-valued mapping $\Phi:\Delta \to {2^{{\Delta_1}}}$ is said to be
\begin{itemize}
\item[(i)]   upper semicontinuous (u.s.c.) at ${u_0} \in \Delta$ if,
for any neighborhood $V$ of $\Phi\left( {{u_0}} \right)$, there exists
a neighborhood $U\left( {{u_0}} \right)$ of ${u_0}$ such that for every $u \in U\left( {{u_0}} \right)$, $\Phi\left( u \right) \subseteq  V$.
\item[(ii)] lower semicontinuous (l.s.c.) at ${u_0} \in \Delta$ if, for any
$x \in \Phi\left( {{u_0}} \right)$ and any neighborhood $V$ of
$x$, there exists a neighborhood $U\left( {{u_0}} \right)$ of
${u_0}$ such that for every $u \in U\left( {{u_0}} \right)$, $\Phi\left( u \right) \cap V \ne \emptyset$.
\item[(iii)]  Hausdorff lower semicontinuous (H-l.s.c.) at ${u_0} \in T$ if,
for any neighborhood $V$ of $0 \in {T_1}$, there exists
a neighborhood $U\left( {{u_0}} \right)$ of ${u_0}$ such that for every $u \in U\left( {{u_0}} \right)$, $G\left( {u_0} \right)\subseteq G\left( u \right)  + V$.
\item[(iv)]  convex if, the graph of $\Phi$, i.e., $Graph\left( \Phi  \right): = \left\{ {\left( {x,y} \right) \in \Delta  \times {\Delta _1}:y \in \Phi \left( x \right)} \right\}$ is a convex set in $\Delta  \times {\Delta _1}$.
\item[(v)] rotund if, $Graph\left( \Phi  \right)$ is convex and for any $\left( {{x_1},{y_1}} \right),\left( {{x_2},{y_2}} \right) \in Graph\left( \Phi  \right)$ with ${x_1} \ne {x_2}$, ${y_1} \ne {y_2}$, we have
$$\left( {\left( {{x_1},{y_1}} \right),\left( {{x_2},{y_2}} \right)} \right) \cap {\left( {\partial Graph\left( \Phi  \right)} \right)^c} \ne \emptyset ,$$
where $\left( {\left( {{x_1},{y_1}} \right),\left( {{x_2},{y_2}} \right)} \right) = \left\{ {\lambda \left( {{x_1},{y_1}} \right) + \left( {1 - \lambda } \right)\left( {{x_2},{y_2}} \right): \lambda  \in \left( {0,1} \right)} \right\}.$
\end{itemize}
We say that $\Phi$ is  u.s.c. and l.s.c. on $\Delta$ if it is  u.s.c. and l.s.c. at each point
$u \in \Delta$, respectively.  $\Phi$ is called to be continuous on $\Delta$ if it is both u.s.c. and l.s.c. on $\Delta$.
\end{definition}
\begin{remark} Obviously, if $\Phi:\Delta \to {2^{{\Delta_1}}}$ is convex, then $\Phi \left( x \right)$ is a convex set for any $x \in \Delta $.
\end{remark}

\begin{lemma} \emph{(\cite{Aubin})}\label{l2.1} A set-valued mapping $\Phi:\Delta \to {2^{{\Delta_1}}}$  is l.s.c. at ${u_0} \in \Delta $ if and only if for any sequence $\left\{ {{u_n}} \right\} \subseteq \Delta $
 with ${u_n} \to {u_0}$ and for any ${x_0} \in \Phi\left( {{u_0}} \right)$, there exists ${x_n} \in \Phi\left( {{u_n}} \right)$ such that ${x_n} \to {x_0}$.
\end {lemma}

\begin{lemma} \emph{(\cite{GHTZ})}\label{l2.2} Let $\Phi:\Delta \to {2^{{\Delta_1}}}$ be a set-valued mapping. For any given $u_0\in \Delta$, if $\Phi\left( {{u_0}} \right)$ is compact, then  $\Phi$ is u.s.c. at ${u_0} \in \Delta $ if and only if    for any sequence $\left\{ {{u_n}} \right\} \subseteq \Delta $
 with ${u_n} \to {u_0}$ and for any ${x_n} \in \Phi\left( {{u_n}} \right)$, there exist ${x_0} \in \Phi\left( {{u_0}} \right)$
 and a subsequence $\left\{ {{x_{{n_k}}}} \right\}$ of $\left\{ {{x_n}} \right\}$ such that ${x_{{n_k}}} \to {x_0}$.
\end {lemma}

\begin{lemma} (Kakutani-Fan-Glicksberg Fixed Point Theorem \cite{Fan2, Glicksberg}). Let $K$ be a nonempty compact convex subset of a locally convex  Hausdorff topological vector space $X$ and let $F:K \to {2^K}$ be an u.s.c. set-valued mapping with nonempty compact convex values.
Then there exists ${x_0} \in K$ such that ${x_0} \in F\left( {{x_0}} \right)$.
\end{lemma}

\section{ The main results}

\begin{lemma} Let $X$ be a reflexive Banach space,  $B$ be the closed unit ball of $X$ and  $A$ be a nonempty closed convex subset of $X$. For given $\delta  > 0$, if $a + \delta B \subseteq A + \delta B$, then $a \in A$.
\end{lemma}
\emph{Proof} Suppose on the contrary that $a \notin A$. Since $A$ is closed, one has
$$d\left( {a,A} \right) = \mathop {\inf }\limits_{y \in A} \left\| {a - y} \right\| > 0.$$
Noting that $X$ is  a reflexive Banach space and $A$ is a nonempty closed convex subset of $X$, there exists $\beta  \in A$ such that
\begin{equation}\label{W1} d\left( {a,A} \right) = \mathop {\inf }\limits_{y \in A} \left\| {a - y} \right\| = \left\| {a - \beta } \right\| > 0.
\end{equation}
Let $\lambda  = {{\left\| {a - \beta } \right\|} \over {\delta  + \left\| {a - \beta } \right\|}}$. We choose $h \in X$ such that $a = \lambda h + \left( {1 - \lambda } \right)\beta $. Then $h = {a \over \lambda } + \beta  - {\beta  \over \lambda }$.

We claim that $\left\| {h - y} \right\| \ge \left\| {h - \beta } \right\|$ for any $y \in A$. It follows from (\ref{W1}) that
\begin{equation}\label{W2} \left\| {a - y} \right\| \ge \left\| {a - \beta } \right\|,\;\; \forall y \in A.
\end{equation}
For any $y \in A$, since $y,\beta  \in A$ and $A$ is convex, we have $\lambda y + \left( {1 - \lambda } \right)\beta  \in A$.  By (\ref{W2}), we know that
$\left\| {a - \left( {\lambda y + \left( {1 - \lambda } \right)\beta } \right)} \right\| \ge \left\| {a - \beta } \right\|$ and so
\begin{equation}\label{W3} \left\| {h - y} \right\| \ge \left\| {h - \beta } \right\|,\;\; \forall y \in A.
\end{equation}
On the other hand,
$$\left\| {h - a} \right\| = \left\| {{a \over \lambda } + \beta  - {\beta  \over \lambda } - a} \right\| = \left( {{1 \over \lambda } - 1} \right)\left\| {a - \beta } \right\| = \delta $$
and so $h \in a + \delta B$. Noting that (\ref{W3}) and $$\left\| {h - \beta } \right\| = \left\| {{a \over \lambda } + \beta  - {\beta  \over \lambda } - \beta } \right\| = {1 \over \lambda }\left\| {a - \beta } \right\| = \delta  + \left\| {a - \beta } \right\| > \delta ,$$
 we have $h \notin A + \delta B$ and so $a + \delta B \not\subset A + \delta B$,  which contradicts the assumption that $a + \delta B \subseteq A + \delta B$. This completes the proof.
\hfill$\Box$

\begin{theorem} Let ${\lambda _0} \in \Lambda $ and  $A$ be a nonempty compact convex subset of  a reflexive Banach space $X$.  Assume that $T\left( { \cdot ,{\lambda _0}} \right)$ is rotund and $T\left( { \cdot , \cdot } \right)$ is   continuous on   $A \times \left\{ {{\lambda _0}} \right\}$ with nonempty closed  convex  values.
Then $S\left(  \cdot  \right)$ is l.s.c. at ${\lambda _0}$.
\end{theorem}
  \emph{Proof} Suppose on the contrary that  $S\left(  \cdot  \right)$ is not l.s.c. at ${\lambda _0}$.  Then there
exist a point ${x_0} \in S\left( {{\lambda _0}} \right)$, a neighborhood ${W_0}$ of $0 \in X$ and a sequence $\left\{ {{\lambda _n}} \right\}$ with ${\lambda _n} \to {\lambda _0}$ such that
\begin{equation}\label{F1} \left( {{x_0} + {W_0}} \right) \cap S\left( { \lambda  _n} \right) = \emptyset, \;\;  \forall n \in \mathbb{N}.
\end{equation}
 There are two cases to be considered.

 Case 1. $S\left( { {\lambda _0}} \right)$ is a singleton. For ${x_n} \in S\left( {{\lambda _n}} \right)$, one has
\begin{equation}\label{F2}
{x_n} \in T\left( {{x_n},{\lambda _n}} \right), \quad \forall n \in \mathbb{N}.
 \end{equation}
 Since ${x_n} \in A$ and $A$ is compact, without loss of generality, we can assume that ${x_n} \to \bar{x} \in A$. Noting that $T\left( { \cdot , \cdot } \right)$ is  u.s.c. at $\left( {\bar{x},{\lambda _0}} \right)$, it follows from Lemma \ref{l2.2} and (\ref{F2}) that there exist a point  $x' \in T\left( {\bar{x},{\lambda _0}} \right)$ and a subsequence $\left\{ {{x_{{n_k}}}} \right\}$ of $\left\{ {{x_n}} \right\}$ such that ${x_{{n_k}}} \to {x'}$. By  ${x_n} \to \bar{x}$, we know that $\bar{x} = x'$ and so $\bar{x} = x' \in T\left( {\bar{x},{\lambda _0}} \right)$. This means that $\bar{x} \in S\left( {{\lambda _0}} \right)$. Noting that $S\left( { {\lambda _0}} \right)$ is a singleton, we have $\bar{x} = {x_0}$ and so ${x_n} \to \bar{x} = {x_0}$. Thus, ${x_n} \in {x_0} + {W_0}$ for $n$ large enough. This together with ${x_n} \in S\left( {{\lambda _n}} \right)$ implies that $\left( {{x_0} + {W_0}} \right) \cap S\left( {{\lambda _n}} \right) \not= \emptyset$ for $n$ large enough, which contradicts  (\ref{F1}).

Case 2. $S\left( { {\lambda _0}} \right)$ is not a singleton. Then there exists ${x^*} \in S\left( {{\lambda _0}} \right)$ such that ${x^*} \ne {x_0}$. Since ${x^*},{x_0} \in S\left( {{\lambda _0}} \right)$,   we know that ${x^*} \in T\left( {{x^*},{\lambda _0}} \right)$ and ${x_0} \in T\left( {{x_0},{\lambda _0}} \right)$. Thus,  $\left( {{x^*},{x^*}} \right),\left( {{x_0},{x_0}} \right) \in Graph\left( {T\left( { \cdot ,{\lambda _0}} \right)} \right)$.
Let
$$
x\left( t \right)= t{x^*} + \left( {1 - t} \right){x_0}, \quad \forall t \in [0,1].
$$
Then it is clear that $x\left( t \right) \in A$.
 Since $Graph\left( {T\left( { \cdot ,{\lambda _0}} \right)} \right)$ is rotund, we  can find  ${t_0} \in \left( {0,1} \right)$ such that
\begin{equation}\label{F3} x\left( {{t_0}} \right) \in {x_0} + {W_0}
\end{equation}  and
\begin{equation}\label{F4} \left( {x\left( {{t_0}} \right),x\left( {{t_0}} \right)} \right) \in {\rm{int}}\left( {Graph\left( {T\left( {\cdot,{\lambda _0}} \right)} \right)} \right).
\end{equation}
It follows from (\ref{F4}) that there exists a constant $\delta  > 0$ such that
$$\left( {x\left( {{t_0}} \right),x\left( {{t_0}} \right)} \right) + \delta B \times \delta B \in Graph\left( {T\left( { \cdot ,{\lambda _0}} \right)} \right),$$
where $B$ is the closed unit ball in  $X$. This shows that
\begin{equation}\label{F5} x\left( {{t_0}} \right) + \delta B \subseteq T\left( {x\left( {{t_0}} \right),{\lambda _0}} \right).
\end{equation}
Since $T\left( {x\left( {{t_0}} \right), \cdot } \right)$ is   l.s.c. at ${\lambda _0}$ and $T\left( {x\left( {{t_0}} \right),{\lambda _0}} \right)$ is compact, we can see that $T\left( {x\left( {{t_0}} \right), \cdot } \right)$ is   H-l.s.c. at ${\lambda _0}$. Thus, for $\delta B$, there exists $n_0 \in \mathbb{N}$ large enough such that $T\left( {x\left( {{t_0}} \right),{\lambda _0}} \right) \subseteq T\left( {x\left( {{t_0}} \right),{\lambda _{{n_0}}}} \right) + \delta B$ and so (\ref{F5}) yields  that
\begin{equation}\label{F6} x\left( {{t_0}} \right) + \delta B \subseteq T\left( {x\left( {{t_0}} \right),{\lambda _{{n_0}}}} \right) + \delta B.
\end{equation}
By the convexity and closedness of $T\left( {x\left( {{t_0}} \right),{\lambda _{{n_0}}}} \right)$,  from (\ref{F6}) and Lemma 3.1, we know that $x\left( {{t_0}} \right) \in T\left( {x\left( {{t_0}} \right),{\lambda _{{n_0}}}} \right)$ and so $x\left( {{t_0}} \right) \in S\left( {{\lambda _{{n_0}}}} \right)$. This together with (\ref{F3})   implies that
$x\left( {{t_0}} \right) \in \left( {{x_0} + {W_0}} \right) \cap S\left( {{\lambda _{{n_0}}}} \right)$,  which contradicts (\ref{F1}). This completes the proof.
\hfill$\Box$
\begin{remark} In \cite{Anh1, Anh2, Anh3}, the authors assume that $S\left(  \cdot  \right)$ is l.s.c. on $ \Lambda $, but they do not give any sufficient condition  guaranteeing  that $S\left(  \cdot  \right)$ is l.s.c. on $ \Lambda $. In Theorem 3.1,  we give a sufficient condition for the lower semicontinuity of $S\left(  \cdot  \right)$.
\end{remark}

\begin{theorem} Let $\left( {{\lambda _0},{\varepsilon _0}} \right) \in \Lambda  \times {\mathbb{R}_ + ^0 }$ and  $A$ be a nonempty compact convex subset of a normed vector space $X$.  Assume that $T\left( { \cdot ,{\lambda _0}} \right)$ is convex, $T\left( { \cdot ,{\lambda _0}} \right)$ is   u.s.c on   $A$ with nonempty closed  values and for any $x \in A$, $T\left( {x, \cdot } \right)$ is l.s.c. at ${\lambda _0}$.
Then  $E\left(  \cdot,  \cdot \right)$ is l.s.c. at $\left( {{\lambda _0},{\varepsilon _0}} \right)$.
\end{theorem}
  \emph{Proof} Suppose on the contrary that  $E\left(  \cdot,  \cdot \right)$ is not l.s.c. at $\left( {{\lambda _0},{\varepsilon _0}} \right)$. Then there
exist a point ${x_0} \in E\left( {{\lambda _0},{\varepsilon _0}} \right)$, a neighborhood ${W_0}$ of $0 \in X$ and a sequence $\left\{ {\left( {{\lambda _n},{\varepsilon _n}} \right)} \right\}$ with $\left( {{\lambda _n},{\varepsilon _n}} \right) \to \left( {{\lambda _0},{\varepsilon _0}} \right)$ such that
\begin{equation}\label{E1} \left( {{x_0} + {W_0}} \right) \cap E\left( {{\lambda _n},{\varepsilon _n}} \right) = \emptyset, \;\;  \forall n \in \mathbb{N}.
\end{equation}
 Define a set-valued mapping $Q:{\mathbb{R}_ + } \to {2^{A}}$ by
$$Q\left( \varepsilon  \right) = E\left( {{\lambda _0},\varepsilon } \right) = \left\{ {x \in A:d\left( {x,T\left( {x,{\lambda _0}} \right)} \right) \le \varepsilon } \right\},\;\;\varepsilon  \in {\mathbb{R}_ + }.$$

We claim that ${Q}\left(  \cdot  \right)$ is l.s.c. on ${\mathbb{R}_ + ^0 }$.
  Suppose on the contrary that  there exists ${\varepsilon _0} \in {R}_ + ^0$ such that ${Q}\left(  \cdot  \right)$ is not l.s.c. at ${\varepsilon _0}$. Then there exist a point ${\bar{x}} \in {Q}\left( {{\varepsilon _0}} \right)$, a neighborhood ${U_0}$ of $0 \in X$ and  a sequence $\left\{ {{\varepsilon _n}} \right\}$ with ${\varepsilon _n} \to {\varepsilon _0}$ such that
\begin{equation}\label{E2} \left( {{\bar{x}} + {U_0}} \right) \cap {Q}\left( {{\varepsilon _n}} \right) = \emptyset , \;\;\; \forall n \in \mathbb{N}.
\end{equation}
It is easy to see that, if $0 \le \alpha  \le \beta $, then ${Q}\left( \alpha  \right) \subseteq {Q}\left( \beta  \right)$.
Suppose that ${\varepsilon _0} \le {\varepsilon _n}$. Then ${\bar{x}} \in {Q}\left( {{\varepsilon _0}} \right) \subseteq {Q}\left( {{\varepsilon _n}} \right)$,
which contradicts (\ref{E2}). Thus, we know that ${\varepsilon _0} > {\varepsilon _n}$ for any $n \in \mathbb{N}$.
It follows from Lemma 2.3 and the closedness of ${T\left( {x,{\lambda _0}} \right)}$
  that $$Q\left( 0 \right) = \left\{ {x \in A:d\left( {x,T\left( {x,{\lambda _0}} \right)} \right) = 0} \right\} = \left\{ {x \in A:x \in T\left( {x,{\lambda _0}} \right)} \right\} \ne \emptyset .$$
We choose  $x' \in {Q}\left( 0 \right)$. It follows from ${\varepsilon _n} \to {\varepsilon _0}$ that there exists ${\varepsilon _{n_0}}$ such that
\begin{equation}\label{E3}
\varepsilon'{\bar{x}} + (1-\varepsilon')x' = {\bar{x}} + (1-\varepsilon')\left( {x' - {\bar{x}}} \right) \in {\bar{x}} + {U_0},
\end{equation}
where $\varepsilon'={{{\varepsilon _{{n_0}}}} \over {{\varepsilon _0}}}$.

Now we claim that
$\varepsilon'{\bar{x}} + (1-\varepsilon')x' \in {Q}\left( {{\varepsilon _{{n_0}}}} \right).$
In fact, since   $x' \in {Q}\left( 0 \right)$, we have $x' \in T\left( {x',{\lambda _0}} \right)$ and so
$\left( {x',x'} \right) \in Graph\left( {T\left( { \cdot ,{\lambda _0}} \right)} \right)$.
It follows from ${\bar{x}} \in {Q}\left( {{\varepsilon _0}} \right)$ that $d\left( {{\bar{x}},T\left( {{\bar{x}},{\lambda _0}} \right)} \right) \le {\varepsilon _0}$. Since ${T\left( {{\bar{x}},{\lambda _0}} \right)}$ is compact, there exists $\bar{y} \in T\left( {\bar{x},{\lambda _0}} \right)$ such that
\begin{equation}\label{E5}
d\left( {\bar{x},\bar{y}} \right) = d\left( {\bar{x},T\left( {\bar{x},{\lambda _0}} \right)} \right) \le {\varepsilon _0}.
\end{equation}
Thus, $\bar{y} \in T\left( {\bar{x},{\lambda _0}} \right)$ shows that $\left( {\bar{x},\bar{y}} \right) \in Graph\left( {T\left( { \cdot ,{\lambda _0}} \right)} \right)$. Since $Graph\left( {T\left( { \cdot ,{\lambda _0}} \right)} \right)$ is convex, one has
$$
\left( \varepsilon'\bar{x} + (1-\varepsilon')x', \varepsilon'\bar{y} + (1-\varepsilon')x' \right) \in Graph\left( {T\left( { \cdot ,{\lambda _0}} \right)} \right)
$$
and so
$$
\varepsilon'\bar{y} +  (1-\varepsilon') x' \in T\left( \varepsilon' \bar{x} + (1-\varepsilon')x',\lambda _0 \right).
$$
Thus, it follows from (\ref{E5}) that
\begin{eqnarray*}
& & d\left(\varepsilon'\bar{x} + (1-\varepsilon')x',T\left(  \varepsilon' \bar{x} +  (1-\varepsilon') x', {\lambda _0} \right) \right) \\
&\le& d\left( \varepsilon'\bar{x} + (1-\varepsilon')x',  \varepsilon' \bar{y} +  (1-\varepsilon')x' \right)\\
 &=& \varepsilon' d\left( {\bar{x},\bar{y}} \right) \le {\varepsilon _{{n_0}}}.
 \end{eqnarray*}
This means that $\varepsilon' x +  (1-\varepsilon') x' \in Q\left( {{\varepsilon _{{n_0}}}} \right)$.
It follows from (\ref{E3}) that
$$
\varepsilon' \bar{x} +  (1-\varepsilon')x' \in \left( {\bar{x} + {U_0}} \right) \cap Q\left( {{\varepsilon _{{n_0}}}} \right),
$$
which contradicts (\ref{E2}). Therefore, ${Q}\left(  \cdot  \right)$ is l.s.c. on ${\mathbb{R}_ + ^0 }$.

For the above ${x_0} \in E\left( {{\lambda _0},{\varepsilon _0}} \right) = Q\left( {{\varepsilon _0}} \right)$ and $W_0$, there exists a  neighborhood $U\left( {{\varepsilon _0}} \right)$ of ${\varepsilon _0}$ such that
$$\left( {{x_0} + {W_0}} \right) \cap Q\left( \varepsilon  \right) \ne \emptyset ,\;\; \forall \varepsilon  \in U\left( {{\varepsilon _0}} \right).$$
Choose ${\varepsilon ^*} \in U\left( {{\varepsilon _0}} \right)$ with $0 < {\varepsilon ^*} < {\varepsilon _0}$. Then
$$\left( {{x_0} + {W_0}} \right) \cap E\left( {{u_0},{\varepsilon ^*}} \right) = \left( {{x_0} + {W_0}} \right) \cap Q\left( {{\varepsilon ^*}} \right) \ne \emptyset $$
and so there exists ${x^*} \in {x_0} + {W_0}$ such that
\begin{equation}\label{E6} d\left( {{x^*},T\left( {{x^*},{\lambda _0}} \right)} \right) \le {\varepsilon ^*}. \end{equation}
It follows from ${x^*} \in {x_0} + {W_0}$ and (\ref{E1}) that ${x^*} \notin E\left( {{\lambda _n},{\varepsilon _n}} \right)$ and so
\begin{equation}\label{E6++} d\left( {{x^*},T\left( {{x^*},{\lambda _n}} \right)} \right) > {\varepsilon _n}. \end{equation}
Let $\delta  = {{{\varepsilon _0} - {\varepsilon ^*}} \over 2} > 0$. By (\ref{E6}), we know that there exists ${y^*} \in T\left( {{x^*},{\lambda _0}} \right)$ such that
 \begin{equation}\label{E7} d\left( {{x^*},{y^*}} \right) < {\varepsilon ^*} + \delta .
 \end{equation}
Since $T\left( {{x^*}, \cdot } \right)$ is l.s.c. at ${\lambda _0}$,  by Lemma \ref{l2.1},    there exists ${y_n} \in T\left( {{x^*},{\lambda _n}} \right)$ such that ${y_n} \to {y^*}$ and so $d\left( {{x^*},{y_n}} \right) \to d\left( {{x^*},{y^*}} \right)$. It follows from (\ref{E7}) that
 \begin{equation}\label{E8} d\left( {{x^*},{y_n}} \right) < {\varepsilon ^*} + \delta
 \end{equation}
for $n$ large enough.  On the other hand,  from (\ref{E6++}) and ${\varepsilon _n} \to {\varepsilon _0}$, we have
$$d\left( {{x^*},{y_n}} \right) \ge d\left( {{x^*},T\left( {{x^*},{\lambda _n}} \right)} \right) > {\varepsilon _n} > {\varepsilon _0} - \delta  = {\varepsilon ^*} + \delta $$
for $n$ large enough, which contradicts (\ref{E8}).  This completes the proof.
\hfill$\Box$

Next, we give an example to illustrate Theorems 3.1  and 3.2.
\begin{example}
 Let  $A = \left[ {0,2} \right]$ and $\Lambda  = \mathbb{R}$.   Let $T:A \times \Lambda  \to {2^A}$ be defined as follows:
$$T\left( {x,\lambda } \right) = \left\{ {y \in \mathbb{R}:\left( {{1 \over 3}\cos \lambda  + {2 \over 3}} \right)\left( {1 - \sqrt {{x^2} + 2x} } \right) \le y \le \left( {{1 \over 3}\cos \lambda  + {2 \over 3}} \right)\left( {1 + \sqrt {{x^2} + 2x} } \right)} \right\}.$$
Let ${\lambda _0} = 0$ and ${\varepsilon _0} >0$. Clearly,
$$Graph\left( {T\left( { \cdot ,{\lambda _0}} \right)} \right) = \left\{ {\left( {x,y} \right) \in {{\mathbb{R}}^2}:{{\left( {x - 1} \right)}^2} + {{\left( {y - 1} \right)}^2} \le 1} \right\}$$
is rotund. It is easy to check that all conditions of Theorems 3.1  and 3.2 are satisfied.  Thus, Theorem 3.1 shows that $S\left(  \cdot  \right)$ is l.s.c. at ${\lambda _0}$ and Theorem 3.2 implies that  $E\left(  \cdot,  \cdot \right)$ is l.s.c. at $\left( {{\lambda _0},{\varepsilon _0}} \right)$.
\end{example}

\section{Applications}

\subsection{Existence of solutions for a class of generalized Stackelberg equilibrium problems}

Following Nagy \cite{Nagy}, assume that  $f_1,f_2: \mathbb{R}^N\times \mathbb{R}^N \to \mathbb{R}$ are the payoff/loss functions for two players, and $K_1,K_2\subset \mathbb{R}^N$ are their strategy sets. It is well known that the framework of Stackelberg equilibrium problem can be modelled  by the following bi-level mathematical programming problem:
\begin{eqnarray*}
&&\min f_1(x,y)\\
&& {\rm s.t.} \quad y\in R_{SE}(x), \; x\in K_1,
\end{eqnarray*}
where $R_{SE}(x)$ denotes the Stackelberg equilibrium response set given  by
$$
R_{SE}(x)=\{y\in K_2: f_2(x,v)-f_2(x,y)\ge 0,\; \forall v\in K_2\}.
$$
Applying the variational inequality theory and the fixed point theorem,  Nagy \cite{Nagy} studied the existence and location of Stackelberg equilibrium problem under the assumptions that $f_1$ and $f_2$ are both smooth functions.  Moreover, Han and Huang \cite{HH} showed the existence of solutions for the Stackelberg
equilibrium problem without the smoothness by employing the lower semicontinuity of the set-valued mapping $R_{SE}(x)$.

Let $K_1,K_2\subset \mathbb{R}^N$ be two subsets and $f : \mathbb{R}^N\times \mathbb{R}^N \to \mathbb{R}$ be a function.  Assume that
$T:{K_2} \times {K_1} \to {2^{{K_2}}}$  is a set-valued mapping. We consider the following  generalized Stackelberg equilibrium problem:
\begin{eqnarray*}
&&\min f(x,y)\\
&& {\rm s.t.} \quad y\in R'_{SE}(x), \; x\in K_1,
\end{eqnarray*}
where $R'_{SE}(x)$ is the Stackelberg equilibrium response set defined by
$${R'_{SE}}(x) = \left\{ {y \in {K_2}:y \in T\left( {y,x} \right)} \right\}.$$

It is well known that, when $f_2$ is smooth or subdifferential,  $y\in R_{SE}(x)$ if and only if $y$ is a fixed point of the single-valued or set-valued mapping. Therefore, the generalized Stackelberg equilibrium problem can be regarded as a generalization of the Stackelberg equilibrium problem considered by Nagy \cite{Nagy}.

Now we are going to give an existence result concerned with the solutions of the generalized Stackelberg equilibrium problem.  To this end, we assume that $K_2$ is a nonempty compact convex subset of $\mathbb{R}^N$, $T\left( { \cdot ,{x}} \right)$ is rotund for any $x \in {K_1}$, and $T\left( { \cdot , \cdot } \right)$ is   continuous on   $K_2 \times K_1$ with nonempty closed  values.  By Lemma 2.3, we can   see that $R'_{SE}(x)$ is nonempty for any $x \in {K_1}$. Since $T\left( { \cdot ,{x}} \right)$ is convex, it is easy to see that $R'_{SE}(x)$ is convex for any $x \in {K_1}$.

 We claim that $R'_{SE}(x)$ is closed for any $x \in {K_1}$. In fact, let $\left\{ {{y_n}} \right\} \subseteq {R'_{SE}}\left( x \right)$ with ${y_n} \to {y_0}$. Then ${y_0} \in {K_2}$ and ${y_n} \in T\left( {{y_n},x} \right)$.
Noting that $T\left( { \cdot , x } \right)$ is  u.s.c. at $y_0$, it follows from Lemma \ref{l2.2}  that there exist a point $y' \in T\left( {{y_0},x} \right)$ and a subsequence $\left\{ {{y_{{n_k}}}} \right\}$ of $\left\{ {{y_n}} \right\}$ such that ${y_{{n_k}}} \to {y'}$. By  ${y_n} \to {y_0}$, we have ${y_0} = y'$ and so ${y_0} = y' \in T\left( {{y_0},x} \right)$. This means that ${y_0} \in {R'_{SE}}\left( x \right)$ and so $R'_{SE}(x)$ is closed.

Therefore, it follows from Theorem 3.1 that $R'_{SE}:K_1\to 2^{K_2}$ is lower semi-continuous on $K_1$. By Michael's continuous selection theorem (see, for example, Theorem 16.1 of \cite{G99}),  there exists a continuous selection $r(x)\in R'_{SE}(x)$ for all $x\in K_1$.

Moreover, we assume that $K_1$ is a nonempty compact  subset of $\mathbb{R}^N$ and $f$ is continuous on $K_1\times K_2$. Then it is easy to see that $f (x,r(x))$ is continuous on $K_1$ and so there exists a point $x^*\in K_1$ such that
$$
f (x^*,r(x^*))=\min_{x\in K_1}f (x,r(x)).
$$
Let $y^*=r(x^*)$. Then we know that $(x^*,y^*)$ is a solution of the generalized Stackelberg equilibrium problem.

\subsection{Lower semicontinuity of the solution mapping for  the parametric vector quasiequilibrium problems}

Let $\Lambda$, $\Omega $ and  $Y$   be three normed vector spaces.
 Let $D$ be a nonempty  subset of  a reflexive Banach space $X$.
Let $K:D \times \Lambda  \to {2^D}$ and $F:D \times D \times \Omega  \to {2^Y}$
 be two set-valued mappings. Let $C \subseteq Y$ be closed with ${\mathop{\rm int}} C \ne \emptyset $. For $\left( {u,\lambda } \right) \in \Omega  \times \Lambda $,
 we consider the following parametric vector quasiequilibrium problems:

\begin{center}(QEP) \quad  finding   ${x_0} \in {\rm{cl}}K\left( {{x_0},\lambda } \right)$ such that \quad $F\left( {{x_0},y,u} \right) \cap \left( {Y\backslash  - {\mathop{\rm int}} C} \right) \ne \emptyset $, \quad ${\rm{y}} \in K\left( {{x_0},\lambda } \right)$;
 \end{center}
\begin{center}(SQEP)  \quad  finding   ${x_0} \in {\rm{cl}}K\left( {{x_0},\lambda } \right)$ such that  \quad $F\left( {{x_0},y,u} \right) \subseteq Y\backslash  - {\mathop{\rm int}} C$, \quad ${\rm{y}} \in K\left( {{x_0},\lambda } \right)$.
 \end{center}

For $\left( {u,\lambda } \right) \in \Omega  \times \Lambda $, let ${M_1}\left( {u,\lambda } \right)$
 denote the set of all  solutions of (QEP), i.e.
$${M_1}\left( {u,\lambda } \right) = \left\{ {x \in {\rm{cl}}K\left( {x,\lambda } \right):F\left( {x,y,u} \right) \cap \left( {Y\backslash  - {\mathop{\rm int}} C} \right) \ne \emptyset ,\;\forall y \in K\left( {x,\lambda } \right)} \right\},$$
and let ${M_2}\left( {u,\lambda } \right)$
 denote the set of all  solutions of (SQEP), i.e.
$${M_2}\left( {u,\lambda } \right) = \left\{ {x \in {\rm{cl}}K\left( {x,\lambda } \right):F\left( {x,y,u} \right) \subseteq Y\backslash  - {\mathop{\rm int}} C,\;\forall y \in K\left( {x,\lambda } \right)} \right\}.$$
For $\lambda  \in \Lambda $, let $H\left( \lambda  \right): = \left\{ {x \in X:x \in {\rm{cl}}K\left( {x,\lambda } \right)} \right\}$. We always assume that ${M_1}\left( {u,\lambda } \right) \ne \emptyset $ and ${M_2}\left( {u,\lambda } \right) \ne \emptyset $ for all $\lambda $ in a neighborhood of ${\lambda _0} \in \Lambda $ and for all $u$ in a neighborhood of ${u _0} \in \Omega$.

Form Theorem 2.1 of \cite{Anh1} and Theorem 3.1, we can get the following theorem.

\begin{theorem}
Let $\left( {{u_0},{\lambda _0}} \right) \in \Omega  \times \Lambda $ and  $D$ be  compact and convex.  Assume that
\begin{itemize}
\item[(i)] $K\left(  \cdot , {\lambda _0} \right)$ is rotund and $K\left(  \cdot , \cdot \right)$ is  continuous on  $D \times \left\{ {{\lambda _0}} \right\}$ with nonempty closed convex values;
\item[(ii)] $F\left( { \cdot , \cdot , \cdot } \right)$ is l.s.c. on $D \times D \times \left\{ {{u_0}} \right\}$;
 \item[(iii)] for any $x \in {M_1}\left( {{u_0}, {\lambda _0}} \right)$ and any $y \in K\left( {x,{\lambda _0}} \right)$, $F\left( {x,y,{u_0}} \right) \cap \left( {Y\backslash  - C} \right) \ne \emptyset $.
\end{itemize}
 Then ${M_1}\left( {  \cdot , \cdot } \right)$ is l.s.c. at $\left( {{u_0},{\lambda _0}} \right)$.
\end {theorem}

\begin{definition}  (\cite{Anh1}) Let $X$ be a  Hausdorff topological space, $Y$ be a topological vector space and $C \subseteq Y$ with ${\mathop{\rm int}} C \ne \emptyset $. A set-valued mapping $\Phi:X \to {2^Y}$ is said to be  have
\begin{itemize}
\item[(i)] the $C$-inclusion property at $x_0$ if, for any ${x_\alpha } \to {x_0}$,
$$\Phi \left( {{x_0}} \right) \cap \left( {Y\backslash  - {\mathop{\rm int}} C} \right) \ne \emptyset  \Rightarrow \exists {\bar{\alpha}} ,\Phi \left( {{x_{\bar{\alpha}} }} \right) \cap \left( {Y\backslash  - {\mathop{\rm int}} C} \right) \ne \emptyset .$$
\item[(ii)]  the strict $C$-inclusion property at $x_0$ if, for any ${x_\alpha } \to {x_0}$,
$$\Phi \left( {{x_0}} \right) \subseteq Y\backslash  - {\mathop{\rm int}} C \Rightarrow \exists {\bar{\alpha}} ,\Phi \left( {{x_{\bar{\alpha}} }} \right) \subseteq Y\backslash  - {\mathop{\rm int}} C.$$
\end{itemize}
\end{definition}

Form Theorem 2.2 of \cite{Anh1} and Theorem 3.1, we can get the following theorem.

\begin{theorem}
Let $\left( {{u_0},{\lambda _0}} \right) \in \Omega  \times \Lambda $ and  $D$ be  compact and convex.  Assume that
\begin{itemize}
\item[(i)] $K\left(  \cdot , {\lambda _0} \right)$ is rotund and $K\left(  \cdot , \cdot \right)$ is  continuous on  $D \times \left\{ {{\lambda _0}} \right\}$ with nonempty closed convex values;
\item[(ii)] $F\left( { \cdot , \cdot , \cdot } \right)$ has the $C$-inclusion property on $D \times D \times \left\{ {{u_0}} \right\}$.
\end{itemize}
 Then ${M_1}\left( {  \cdot , \cdot } \right)$ is l.s.c. at $\left( {{u_0},{\lambda _0}} \right)$.
\end {theorem}

Form Theorem 2.3 of \cite{Anh1} and Theorem 3.1, we can get the following theorem.

\begin{theorem}
Let $\left( {{u_0},{\lambda _0}} \right) \in \Omega  \times \Lambda $ and  $D$ be  compact and convex.  Assume that
\begin{itemize}
\item[(i)] $K\left(  \cdot , {\lambda _0} \right)$ is rotund and $K\left(  \cdot , \cdot \right)$ is  continuous on  $D \times \left\{ {{\lambda _0}} \right\}$ with nonempty closed convex values;
\item[(ii)] $F\left( { \cdot , \cdot , \cdot } \right)$ is u.s.c. on $D \times D \times \left\{ {{u_0}} \right\}$;
 \item[(iii)] for any $x \in {M_2}\left( {{u_0}, {\lambda _0}} \right)$ and any $y \in K\left( {x,{\lambda _0}} \right)$, $F\left( {x,y,{u_0}} \right) \subseteq Y\backslash  - C$.
\end{itemize}
 Then ${M_2}\left( {  \cdot , \cdot } \right)$ is l.s.c. at $\left( {{u_0},{\lambda _0}} \right)$.
\end {theorem}

Form Theorem 2.4 of \cite{Anh1} and Theorem 3.1, we can get the following theorem.

\begin{theorem}
Let $\left( {{u_0},{\lambda _0}} \right) \in \Omega  \times \Lambda $ and  $D$ be  compact and convex.  Assume that
\begin{itemize}
\item[(i)] $K\left(  \cdot , {\lambda _0} \right)$ is rotund and $K\left(  \cdot , \cdot \right)$ is  continuous on  $D \times \left\{ {{\lambda _0}} \right\}$ with nonempty closed convex values;
\item[(ii)] $F\left( { \cdot , \cdot , \cdot } \right)$ has the strict $C$-inclusion property on $D \times D \times \left\{ {{u_0}} \right\}$.
\end{itemize}
 Then ${M_2}\left( {  \cdot , \cdot } \right)$ is l.s.c. at $\left( {{u_0},{\lambda _0}} \right)$.
\end {theorem}

\section*{Acknowledgements}
The authors are grateful to the editor and the referees for their valuable comments and suggestions.

\end{document}